\documentclass[11pt,a4paper]{article}

\usepackage{latexsym,amssymb, amsthm}
\usepackage[centertags]{amsmath}
\usepackage{epsfig,pifont}
\usepackage[inactive]{srcltx}  

\usepackage{graphics,helvet,times,mathptm,epic,eepic}

\usepackage{color}

\usepackage{newlfont}

\usepackage{amsfonts,bbm}

\newenvironment{keywords}{ \noindent {\small\bf Key Words}:}{ }

\def\bd{\begin{description}}
\def\ed{\end{description}}

\def\beq{\begin{equation}}
\def\eeq{\end{equation}}
\def\bea{\begin{eqnarray}}
\def\eea{\end{eqnarray}}
\def\beas{\begin{eqnarray*}}
\def\eeas{\end{eqnarray*}}

\def\G1{\hbox{$\displaystyle{\mbox{\ding{172}}}$}}

\newtheorem{lemma}{Lemma}[section]
\newtheorem{theorem}{Theorem}[section]
\newtheorem{corollary}{Corollary}[section]

\theoremstyle{remark}

\def\ligne#1{\hbox to \hsize{#1}}
\begin{document}

\title{\textbf{\textsc{On existence of infinite primes and infinite twin primes}}
}

\newcommand{\nms}{\normalsize}
\author{  {  \bf Maurice Margenstern\footnote{Maurice Margenstern, Ph.D.,   is Professor
Emeritus at the Universit\'e de Lorraine,
 e-mail: {\tt
margenstern@gmail.com }} \hspace{2mm}and Yaroslav D.
Sergeyev\footnote{Yaroslav D. Sergeyev, Ph.D., D.Sc., D.H.C., is
Distinguished Professor at the University of Calabria, Rende, Italy.
 He is also  Professor (a part-time contract) at the  Lobatchevsky State University,
  Nizhni Novgorod, Russia and Affiliated Researcher at the Institute of High Performance
  Computing and Networking of the National Research Council of Italy, e-mail: {\tt  yaro@si.dimes.unical.it}}
  }\,\,\thanks{Corresponding author.} }

\date{}

\maketitle

\begin{abstract}
The twin primes conjecture is a very old problem. Tacitly it is
supposed that the  primes it deals with are finite. In the present
paper we consider three problems that are not related to finite
primes but deal with infinite integers. The main tool of our
investigation is a  numeral system proposed recently that allows one
to express various infinities and infinitesimals easily and by a
finite number of symbols. The problems under consideration are the
following and for all of them we give affirmative answers: (i) do
infinite primes exist? (ii) do infinite twin primes exist? (ii) is
the set of infinite twin primes infinite? Examples of these three
kinds of objects are given.
 \end{abstract}

\begin{keywords}
Infinite primes, infinite twin primes, infinite  sets.
 \end{keywords}

\section{Introduction}
\label{s1}

   The research on twin primes number has given rise to an important number of works
 (see, for example,
\cite{Guy} and references given therein). In this paper, we study
the problem of the existence of infinite primes and infinite twin
primes. So, our paper brings in no new light on the traditional twin
primes conjecture. Let us  mention a few works which can be
considered as precursors in some sense and might be thought of an
exotic character by pure mathematicians. As an example of such an
exotic issue, we can quote the existence of other natural families
of numbers whose distribution is alike that of primes. There is an
example of such a family for which an analogue of the twin primes
can be formulated and was indeed proved, see~\cite{mmpractical}.
Now, what we can consider as precursors for us are more connected
with logical problems. The first paper in this direction
is~\cite{kemeny}, where the possibility of infinite prime numbers of
the form we consider in this paper was investigated without reaching
definite results. The second paper is~\cite{penzin} where the author
constructs a model of a theory in which it is possible to prove the
twin prime conjecture. There is an important difference in this
paper with our result in this sense that the theory considered by
this author assumes a weak induction axiom while here, we have no
induction at all on infinite integers. The form of the twin primes
in that paper has some similarity with that of our infinite twin
primes, although the paper appeals to higher results in algebra and
analysis which is not at all the case of our work.

A   numeral system introduced recently in
\cite{Sergeyev,informatica}  for performing computations with
infinities and infinitesimals is used here to study the problem of
the existence of infinite primes and infinite twin primes.  It
should be mentioned that this computational methodology   is not
related to non-standard analysis of Robinson and   has a strong
applied character. In fact, the Infinity Computer working
numerically with a variety
 of
infinite and infinitesimal numbers has been introduced  (see the
 patent \cite{Sergeyev_patent}).

  In order to see the place
of the new approach in the historical panorama of ideas dealing with
infinite and infinitesimal, see
\cite{Kauffman,Lolli,Lolli_2,MM_bijection,Sorbi,Dif_Calculus,
first,Sergeyev_Garro}. In particular, connections of the new
approach with   bijections is studied in \cite{MM_bijection} and
metamathematical investigations on the theory can be found in
\cite{Lolli_2}. Among the applications where the new approach has
been successfully used we can mention the following: percolation and
biological processes (see \cite{Iudin,Iudin_2,DeBartolo,Biology}),
hyperbolic geometry (see \cite{Margenstern,Margenstern_3}),
numerical differentiation and optimization (see
\cite{DeLeone,Num_dif,Zilinskas}), infinite series (see
\cite{Kanovei,Dif_Calculus,Riemann,Zhigljavsky}), the first Hilbert
problem, Turing machines, and lexicographic ordering  (see
\cite{first,medals,Sergeyev_Garro,Sergeyev_Garro_2}), cellular
automata (see \cite{DAlotto,DAlotto_3,DAlotto_2}),  ordinary
differential equations (see \cite{ODE}), etc.

The rest of the paper has the following structure. The next Section
contains a brief informal description of the  numeral system
allowing one to express different infinities and infinitesimals in a
unique framework (see \cite{informatica,Lagrange,chapter} for a
detailed discussion). Section~3 presents main results of the paper.

\section{Infinities and infinitesimals expressed in grossone-based numerals }
\label{s2}

Let us consider a study published in \textit{Science}  (see
\cite{Gordon}) where there is a description of a primitive tribe
living in Amazonia - Pirah\~{a} - that uses a very simple numeral
system\footnote{ We remind that \textit{numeral}  is a symbol or a
group of symbols that represents a \textit{number}. The difference
between numerals and numbers is the same as the difference between
words and the things they refer to. A \textit{number} is a concept
that a \textit{numeral} expresses. The same number can be
represented by different numerals. For example, the symbols `3',
`three', and `III' are different numerals, but they all represent
the same number.} for counting: one, two, many.

For Pirah\~{a}, all quantities larger than two are just `many' and
such operations as 2+2 and 2+1 give the same result, i.e., `many'.
Using their weak numeral system Pirah\~{a} are not able to see
numbers 3, 4, etc., to execute arithmetical operations with them,
and, in general, to say anything about these numbers because in
their language there are neither words nor concepts for that.
Moreover, the weakness of their numeral system leads to such results
as
\[
\mbox{`many'}+ 1= \mbox{`many'},   \hspace{1cm}    \mbox{`many'} + 2
= \mbox{`many'},
\]
which are very familiar to us  in the context of views on infinity
used in the traditional calculus
\[
\infty + 1= \infty,    \hspace{1cm}    \infty + 2 = \infty.
\]
 This analogy  advices that
our difficulty in working with infinity is not connected to the
nature of infinity but is just a result of inadequate numeral
systems used to express numbers. In fact,  numeral systems strongly
influence our capabilities to describe physical and mathematical
objects. For instance, Roman numeral system has no numeral to
express 0. As a consequence, the expression III-X in this numeral
system is an indeterminate form. Moreover, any assertion regarding
negative numbers and zero cannot be formulated using Roman numerals
because there are no symbols corresponding to these concepts in this
concrete numeral system.

The  numeral system proposed in \cite{informatica,Lagrange,chapter}
is based on an infinite unit of measure expressed by the numeral \G1
called \textit{grossone} and introduced as the number of elements of
the set $\mathbb{N}$ of natural numbers (a clear difference with
non-standard analysis can be seen immediately since non-standard
infinite numbers are not connected to concrete infinite sets and do
not belong to $\mathbb{N}$). Other symbols dealing with infinities
and infinitesimals ($\infty$, Cantor's $\omega$, $\aleph_0,
\aleph_1, ...$, etc.)  are not used together with \G1. Similarly,
when the positional numeral system and the numeral 0 expressing zero
had been introduced, symbols V, X, and other symbols from the Roman
numeral system had not been involved.

Notice that people very often do  not pay a great attention to the
distinction between numbers and numerals (in this occasion it is
necessary to recall constructivists who studied this issue), many
theories dealing with infinite and infinitesimal quantities have a
symbolic (not numerical) character. For instance, many versions of
non-standard analysis are symbolic, since they have no numeral
systems to express their numbers by a finite number of symbols (the
finiteness of the number of symbols is necessary for organizing
numerical computations). Namely, if we consider a finite $n$ than it
can be taken $n=7$, or $n=108$ or any other numeral used to express
finite quantities and consisting of a finite number of symbols. In
contrast, if we consider a non-standard infinite $m$ then it is not
clear which numerals can be used to assign a concrete value to $m$.
One of the important differences between the new approach and
non-standard analysis consists of the fact that the new numeral
system allows us to assign concrete values to infinities (and
infinitesimals) as it happens with finite values. In fact, we can
assign $m=\G1$, $m=3\G1-2$ or to use any other infinite numeral
involving grossone to give a numerical value to $m$ (see
\cite{informatica,Lagrange,chapter} for a detailed discussion).

 The
numeral \G1 allows one to construct different numerals expressing
different infinities and infinitesimals and to execute numerical
computations with all of them. As a result,  in  occasions requiring
infinities and infinitesimals indeterminate forms and various kind
of divergence are not present when one works with any (finite,
infinite, or infinitesimal) numbers expressible in the new numeral
system and it becomes possible to execute arithmetical operations
with a variety of different infinities and infinitesimals. For
example, for $\mbox{\ding{172}}$ and $\mbox{\ding{172}}^{3.1}$ (that
are examples of infinities) and $\mbox{\ding{172}}^{-1}$ and
$\mbox{\ding{172}}^{-3.1}$ (that are examples of infinitesimals) it
follows
 \beq
 0 \cdot \mbox{\ding{172}} =
\mbox{\ding{172}} \cdot 0 = 0, \hspace{3mm}
\mbox{\ding{172}}-\mbox{\ding{172}}= 0,\hspace{3mm}
\frac{\mbox{\ding{172}}}{\mbox{\ding{172}}}=1, \hspace{3mm}
\mbox{\ding{172}}^0=1, \hspace{3mm} 1^{\mbox{\tiny{\ding{172}}}}=1,
\hspace{3mm} 0^{\mbox{\tiny{\ding{172}}}}=0,
 \label{3.2.1}
       \eeq
\[
 0 \cdot \mbox{\ding{172}}^{-1} =
\mbox{\ding{172}}^{-1} \cdot 0 = 0, \hspace{5mm}
\mbox{\ding{172}}^{3.1} > \mbox{\ding{172}}^{1} > 1 >
\mbox{\ding{172}}^{-1} > \mbox{\ding{172}}^{-3.1} > 0,
\]
\[
\mbox{\ding{172}}^{-1}-\mbox{\ding{172}}^{-1}=
0,\hspace{5mm}\frac{\mbox{\ding{172}}^{-1}}{\mbox{\ding{172}}^{-1}}=1,
\hspace{2mm}
\frac{5+\mbox{\ding{172}}^{-3.1}}{\mbox{\ding{172}}^{-3.1}}=5\mbox{\ding{172}}^{3.1}+1,
\hspace{2mm} (\mbox{\ding{172}}^{-1})^0=1,
       \]
       \[
 \mbox{\ding{172}} \cdot
\mbox{\ding{172}}^{-1} =1,\hspace{2mm} \mbox{\ding{172}} \cdot
\mbox{\ding{172}}^{-3.1}
=\mbox{\ding{172}}^{-2.1},\hspace{2mm}\frac{\mbox{\ding{172}}^{3.1}+4\G1}{\mbox{\ding{172}}}=\G1^{2.1}+4,
       \]
       \[
 \frac{\mbox{\ding{172}}^{3.1}}{\mbox{\ding{172}}^{-3.1}}=\G1^{6.2},\hspace{3mm}       (\mbox{\ding{172}}^{3.1})^0=1, \hspace{3mm}
\mbox{\ding{172}}^{3.1} \cdot \mbox{\ding{172}}^{-1} =\G1^{2.1},
\hspace{3mm} \mbox{\ding{172}}^{3.1} \cdot \mbox{\ding{172}}^{-3.1}
=1.
       \]

  It follows from (\ref{3.2.1}) that
$\mbox{\ding{172}}^0=1$, therefore, a finite number $a$ can be
represented in the new numeral system simply as
$a\mbox{\ding{172}}^0=a$, where the numeral $a$ itself can be
written down by any convenient numeral system used to express finite
numbers. The simplest infinitesimal numbers are represented by
numerals having only negative finite   powers of \G1 (e.g.,
$50.1\mbox{\ding{172}}^{-10.2}\small{+}
16.38\mbox{\ding{172}}^{-20.3}$, see also examples above). Notice
that all infinitesimals are not equal to zero. In particular,
$\frac{1}{\mbox{\ding{172}}}>0$ because it is a result of division
of two positive numbers. We shall not speak more about
infinitesimals here since they are not used in the present paper.

It should be mentioned that in certain cases  \G1-based numerals
allow us to execute a finer analysis of infinite objects than
traditional tools allow us to do. For instance,  it becomes possible
to measure certain infinite sets and   to see, e.g., that the sets
of even and odd numbers have $\G1/2$ elements each. The set
$\mathbb{Z}$ of integers has $2\G1+1$ elements (\G1 positive
elements, \G1 negative elements, and zero). Within the countable
sets and sets having cardinality of the continuum (see
\cite{Lolli,first,Lagrange}) it becomes possible to distinguish
infinite sets having different number of elements expressible in the
numeral system using grossone and to see that, for instance,
 \[
 \frac{\mbox{\ding{172}}}{2} < \mbox{\ding{172}}-1 < \mbox{\ding{172}} < \mbox{\ding{172}}+1 < 2\mbox{\ding{172}}+1 <
 2\mbox{\ding{172}}^2-1 < 2\mbox{\ding{172}}^2    <
 2\mbox{\ding{172}}^2+1  <
\]
 \[
2\mbox{\ding{172}}^2+2  <  2^{\mbox{\ding{172}}}-1 <
2^{\mbox{\ding{172}}} < 2^{\mbox{\ding{172}}}+1 <
10^{\mbox{\ding{172}}} <
  \mbox{\ding{172}}^{\mbox{\ding{172}}}-1 <
  \mbox{\ding{172}}^{\mbox{\ding{172}}} <
  \mbox{\ding{172}}^{\mbox{\ding{172}}}+1.
  \]
It is important to stress that the new approach does not contradict
Cantor's results. The situation is similar to what happens when one
uses a   microscope with two different lenses: the first of them is
weak and allows one to see the object of the observation as two dots
and another lens is stronger and allows the observer to see instead
of the first dot 10 dots and instead of the second dot 32 dots. Both
lenses give two correct answers having different accuracies.
Analogously, both approaches, Cantor's and the new one, give correct
answers but the accuracy of the answers is different. Cantor's tools
say that the sets of even, odd, natural, and integer numbers have
the same cardinality $\aleph_0$. This answer is correct with the
precision that cardinal numbers have. However, the fact that they
all have the same cardinality  can be viewed also as the accuracy of
the used instrument  is too low to see that these sets have
different numbers of elements. The new numeral system allows us to
see these differences among sets having cardinality $\aleph_0$ and
among sets having   cardinality of the continuum, as well (see
\cite{Lolli,MM_bijection,first,Lagrange,chapter} for a detailed
discussion including also the one-to-one correspondence issues).

We conclude this brief informal introduction by mentioning that
properties of grossone are described by the Infinite Unit Axiom (see
\cite{first,Lagrange,chapter}) that is added to axioms for real
numbers. In the context of the present paper two issues postulated
by the Axiom are important for us: grossone is an infinite number;
(ii) grossone is divisible by any finite integer. Notice that
grossone is not the only number enjoying the latter property. In
fact, zero is also divisible by any finite integer.

\section{Infinite primes and twin primes}
\label{s3}

We need some definitions and conventions to continue our study.
First, in the further consideration we use the following
representation of an infinite number~$c$ where its infinite part is
separated from its finite part: $c = c_1+c_2$. In this separation,
$c_1$ is infinite and is expressed by numerals involving  \G1 and
$c_2$ is   finite and is represented by numerals used to write down
finite numbers. Note that such a representation is not unique. For
instance, the number $k=1.7\G1 -1.5$ can be decomposed as
$k_1=1.7\G1,\,\,k_2=-1.5,$ or as
$\tilde{k}_1=1.7\G1-1,\,\,\tilde{k}_2=-0.5,$ or in some other way.
Clearly, this decomposition becomes unique if we require that the
part $c_1$ does not contain any part expressed by finite numerals
only. In our example with the decomposition of the number $k$ we
have its unique decomposition $k_1=1.7\G1,\,\,k_2=-1.5$.

Then, infinite numbers that do not contain finite parts  are called
\textit{purely infinite}. In other words, this means that in their
unique decomposition they have $c_2=0$. Infinite numbers having
their infinite part $c_1$ including more than one infinite part
represented by different powers of \G1 are called \textit{compound}.
Numbers that are not compound are called \textit{simple}. For
instance, the numbers $\G1 - 3\G1^{\frac{1}{2}}$ and $\G1^2+\G1+3.5$
are both compound; the former is purely infinite while the latter is
not. The numbers $\frac{\G1}{2}$ and $\G1^2+1$  are examples of
simple infinite numbers; again the former is purely infinite while
the latter is not. Finally, if a finite or infinite number $c$ is
the square of an integer $d$, i.e., $c=d^2$, we say hereinafter
simply that $c$ is a square.


\begin{lemma} \label{preliminary2}
There exist purely infinite simple numbers $\lambda$ divisible by
all finite integers.
\end{lemma}

\textbf{Proof.} Due to its definition $\G1$ is such a number
$\lambda$. Then, for any positive finite number~$n$, $\frac{\G1}{n}$
is also an example of such a  number $\lambda$.      In fact, for
any finite   number~$p$ the product $pn$ is finite and, therefore,
it follows $pn \vert \G1$ and, as a consequence, $p\vert
\frac{\G1}{n}$. Analogously, $\G1^2$ and $\frac{\G1^2}{n}$ for any
positive finite number~$n$ are examples of   $\lambda$.
  \hfill $\Box$

\begin{theorem}\label{infprimes}
For all purely infinite simple positive integers~$\lambda$ such that
any finite positive integer divides $\lambda$ it follows that
$\lambda+1$ is a prime number.
\end{theorem}

\textbf{Proof.} Let us consider the number $\lambda+1$, where
$\lambda$ is a purely infinite simple positive integer such $p$
divides~$\lambda$, whatever the finite   number~$p$ is. Let us show
that there are no integers $a$ and $b$ such that $a \cdot b =
\lambda +1$.

Suppose that such integers exist. Then two situations are possible:
(i) $a$ is finite and $b$ is infinite; (ii) both $a$ and $b$ are
infinite. The first situation cannot hold because $\lambda$ is
divisible by any finite   number~$p$. Therefore if $p\vert \lambda
+1$, as also $p\vert \lambda$ then $p\vert 1$ which is impossible.
And so, $\lambda+1$ cannot  have a non trivial finite divisor.

  This fact has its importance also for the case (ii)  where $a$
and $b$ are both infinite. It means they they cannot have finite
divisors. Suppose the opposite, i.e., $a=k \cdot c$ where $k$ is a
finite integer and $c$ is an infinite integer. Then $\lambda-1 = k
\cdot c b$, i.e., we have that $\lambda-1$ is  product of the finite
integer $k$ and the infinite integer $cb$. Since $k$ is finite, due
to assumptions of the lemma $k\vert \lambda$ and, therefore,
$k\not\vert \,\, \lambda-1$. The obtained contradiction proves that
infinite numbers $a$ and $b$  cannot have finite divisors.

Let us look at the case (ii) and consider the unique decomposition
of numbers $a$ and $b$ in the form
 \beq
  a = a_1+a_2, \hspace{15mm} b = b_1+b_2,
\label{twins_1}
 \eeq
where the first parts are purely infinite integers (simple or
compound)   and the second parts are finite integers. Then it should
be
 \[
 a \cdot b = (a_1+a_2)(b_1+b_2)= \lambda + 1,
\]
 \beq
   a_1 b_1+ a_2 b_1+ b_2 a_1 + a_2b_2= \lambda + 1.
\label{twins_3}
 \eeq
\noindent

  Let us denote the left-hand part of (\ref{twins_3})   by $L$ and
the right-hand part of (\ref{twins_3})   by $R$.


\noindent Note that $\vert a_1b_1\vert > \vert a_2b_1\vert$ and that
$\vert a_1b_1\vert > \vert b_2a_1\vert$ as $a_2$ and $b_2$ are
finite integers. Also note that $a_2b_2=0$ is impossible as $L$
would contain no finite part while $R$ does. Consequently,
$\displaystyle{\frac{\vert a_1b_1\vert}{\vert a_2b_1\vert}}$ and
$\displaystyle{\frac{\vert a_1b_1\vert}{\vert b_2a_1\vert}}$ are
both infinite numbers. This means that in $L$, the three terms
$a_1b_1$, $a_2b_1$ and $b_2a_1$ are infinite numbers where one of
them, $a_1b_1$, is of a higher order than the others and they do not
contain a finite part in the sense of the unique decomposition
of~\ref{twins_1}. Consequently, $a_2b_2=1$ and $a_2b_1+b_2a_1=0$
since $a_2b_1+b_2a_1$ is a smaller infinite than $a_1b_1$ and
$\lambda$ is a purely infinite simple number. As $a_2$ and$~b_2$ are
integers, we get $a_2=b_2=1$ or $a_2=b_2=-1$. Possibly changing the
sign of $a_1$ and $b_1$ we may assume that $a_2=b_2=1$. This entails
that \hbox{$a_1+b_1=0$} which gives us $-a_1^2 = \lambda$ which is
impossible as $\lambda$ is positive. \hfill $\Box$

\begin{corollary}\label{infiniprimes}
For all finite positive integers~$n$, $\frac{\G1}{n}+1$ and
$\frac{\G1^2}{n}+1$ are  infinite prime numbers.
\end{corollary}

Thus, we have proved the existence of infinite prime numbers and
have given examples of such numbers expressible in \G1-based
numerals. Let us consider now the problem of the existence of
infinite twin primes.


\begin{lemma}\label{infprimes_no}
For all purely infinite simple positive integers~$\lambda$ such that
any finite positive number divides $\lambda$ the infinite integer
$\lambda-1$ is a prime number if and only if~$\lambda$ is not the
square of an integer.
\end{lemma}

\textbf{Proof.}   Let us repeat  the argument of
Theorem~\ref{infprimes} and find the factors of the number
$\lambda-1$. Namely, we have that
  $ab=\lambda-1$, $a=a_1+a_2$, $b=b_1+b_2$ and as in the proof of
Theorem~\ref{infprimes}  we  obtain that $a_2b_2=-1$ and, therefore,
it follows that $a_1=b_1$ and $a_1^2=\lambda$. If $\lambda$ is not a
square, this is impossible and so, we get that $\lambda-1$ is a
prime number. Now, if $\lambda$ is a square,  we obtain $\lambda =
(a_1-1)(a_1+1)$ where the integer $a_1$ satisfies $a_1 = \sqrt{\lambda}$ and,
therefore, $\lambda-1$ cannot be a prime number. \hfill $\Box$

Let us give an example. The infinite integer $\G1^2$ can be taken as
the number $\lambda$ and it can be easily seen that $\G1^2-1$ is not
prime since $\G1^2-1=(\G1 -1)(\G1+1)$ where $\G1 -1$ and $\G1+1$ are
infinite integers.

\begin{lemma}\label{square_no}
For all purely infinite simple positive integers~$\lambda$ such that
any finite positive number divides $\lambda$ and it is a square it
follows that the infinite number $\frac{\lambda}{p^{2m+1}}$ cannot
be a square for finite $m$ and $p$ where $p$ is a prime number.
\end{lemma}

\textbf{Proof.} Suppose that $\frac{\lambda}{p^{2m+1}}$ is a square.
Then it follows $\frac{\lambda}{p^{2m+1}}= r^2$, where, since
$p^{2m+1}$ is a finite   number, $r$ is an infinite integer. Thus,
we can write
 $$\lambda = p^{2m+1}\cdot r^2 = (p^{m}r)^2 \cdot p.$$
Since $p$ is a finite prime number, it cannot be a square. This
result contradicts the fact that $\lambda$ is the square of an
infinite integer and, therefore, we have proved that
$\frac{\lambda}{p^{2m+1}}$ cannot be a square. \hfill $\Box$

Obviously, the infinite number $\frac{\lambda}{2^{2m+1}}$ is an
example illustrating Lemma~\ref{square_no}.

\begin{theorem}  \label{final}
For all purely infinite simple positive integers~$\lambda$ such that
any finite positive number divides $\lambda$ and $\lambda$ is a
square it follows that the numbers $\frac{\lambda}{p^{2m+1}}-1$ and
$\frac{\lambda}{p^{2m+1}}+1$ are infinite twin primes for all
positive finite numbers~$m$  and $p$ where $p$ is a prime number.
\end{theorem}

\textbf{Proof.}   By repeating  the argument of
Theorem~\ref{infprimes} we find that the number
$\frac{\lambda}{p^{2m+1}}+1$ is prime. It follows from
Lemma~\ref{square_no} that $\frac{\lambda}{p^{2m+1}}$ is not a
square. Thus, due to Lemma~\ref{infprimes_no} it follows that
$\frac{\lambda}{p^{2m+1}}-1$ is prime. \hfill $\Box$

Numbers $\frac{\G1^2}{2^{2m+1}}-1$ and $\frac{\G1^2}{2^{2m+1}}+1$
are examples of infinite twin primes for all positive finite
numbers~$m$.

Notice that in   Theorems~\ref{infprimes} and~\ref{final}  we need
only to assume that $\lambda$ is divisible by any finite prime
number. In the following theorem the assumption that any finite
positive integer divides $\lambda$ becomes essential.

\begin{theorem} \label{final_inf}
If:
 \bd
  \item (i) $\lambda$ is an infinite simple positive integer
such that any finite positive integer divides $\lambda$;
 \item (ii) $\lambda$
is a square;
 \item (iii) $m$ is a positive integer;
 \item (iv) $p$ is a finite prime number;
  \ed
   then:
  \bd
  \item (i)
   the sets $A(p)$ have
infinitely many elements, where
 \beq
 A(p)= \{ x: x= \frac{\lambda}{p^{2m+1}},   \,\, p^{2m+1 } \vert \lambda
 \};
\label{twins_4}
 \eeq

  \item (ii)
 all numbers $x \in A(p)$ are infinite integers;
   \item (iii)
   the number of elements of the set $A(p)$ is
 \beq
   M(p)= \max \{ m: y \cdot p^{2m+1} = \lambda \},
   \label{twins_4.1}
 \eeq
   where $p \, \vert  \, y$ and  $p^2 \nmid \,\,  y$.
   \ed
\end{theorem}

\textbf{Proof.} Suppose that the set $A=A(p)$ has $M=M(p)$ elements
and $M$ is finite. Since numbers $p^{2m+1 }$ are strictly
increasing, $p^{2M+1}$ is the largest element in the set. Let us
consider the number $p^{2M+2}$. Since $M$ is finite, it follows that
if \hbox{$p^{2M+2}\vert \lambda$}, it should belong to~$A$. However,
$M+1> M$, thus it cannot belong to $A$. This contradiction concludes
the proof of the first assertion of the theorem.

Let us prove now that all the elements of the set $A$ are infinite
integers. If $m$ in (\ref{twins_4}) is a finite integer then the
respective number $y=\frac{\lambda}{p^{2m+1}}$ is obviously an
infinite integer. Suppose now that $m$ is infinite and $y$ is
finite. Remind that $\lambda$ is divisible by all finite numbers.
Thus, $\frac{\lambda}{p^{2m+1}}$ should be  divisible by all finite
numbers excluding, probably, $p$.  This means that $y$ cannot be
finite since in this case it would be divisible only by a finite
number of integers.

Let us prove  the third assertion of the theorem. The fact $p \,
\vert  \, y$ follows from our supposition that $\lambda$ is a
square. Suppose now that $y = p^2
 y_1$. Then we obtain that
 \[
  y_1 \cdot p^2 \cdot p^{2M+1} = y_1 \cdot     p^{2(M+1)+1} =
  \lambda.
  \]
This contradict the fact that $M$ is the maximal number such that
$p^{2M+1} \vert \lambda$.
 \hfill
 $\Box$

\begin{corollary}\label{final result1}
 Results of Lemma~\ref{square_no} and Theorem~\ref{final} hold for infinite values of~$m$.
\end{corollary}

\begin{corollary}\label{final result2}
 The sets
\beq
 B(p)= \{ x-1,\,\, x+1: x= \frac{\G1^2}{p^{2m+1}},   \,\, p^{2m+1 } \vert \G1^2 \},
\label{twins_5}
 \eeq
of infinite prime numbers have $2M(p)$ elements where the infinite
number $M(p)$ is from (\ref{twins_4.1}).
\end{corollary}

For instance, it follows that the set
 \[
 B(2)= \{ x-1,\,\, x+1: x= \frac{\G1^2}{2^{2m+1}},   \,\, 2^{2m+1 } \vert \G1^2 \},
\]
consists of infinite prime numbers and has infinitely many elements.

 We
conclude the paper with the following rather obvious remark:
substituting $\G1^2$ in (\ref{twins_5}) by $\G1^4,\G1^{16}$ or by
any other  infinite simple positive integer $\lambda$ being  a
square we can generate other infinite sets of infinite prime
numbers.

\bibliographystyle{plain}
\bibliography{MMYS_Twin_primes}
\end{document}